# SDRE Based Attitude Control Using Modified Rodriguez Parameters

R. Ozgur DORUK, Middle East Technical University, Kalkanli, Guzelyurt, Mersin 10, TURKEY

e – mail: rdoruk@metu.edu.tr

*Abstract:*

The purpose of this paper is to present an application of the State Dependent Riccati Equation (SDRE) method to satellite attitude control where the satellite kinematics is modeled by Modified Rodriguez Parameters (MRP). The SDRE methodology is applicable on special forms of nonlinear systems where satellite model is one of the candidates. It is not easy to find an analytical solution from the SDRE. Thus point wise solutions are interpolated with respect to the operating conditions. The point wise solutions are obtained from the MATLAB algorithms which are derived from the positive definite solutions of the SDRE. The global stability analysis is difficult due to the nature of the methodology. The resultant attitude controllers outside the breakpoints (the selected operating conditions for interpolation) are suboptimal. The performance of the designs is examined by simulations on MATLAB – Simulink environment. The simulation results show that, the designed attitude controllers are working satisfactorily even in the presence of inertial uncertainties.

## Introduction

Satellite attitude control gains importance as the technology evolves. However, the nonlinearity of the satellite model makes the process difficult. One solution is the linearization where the nonlinear satellite models are linearized around the origin (zero attitude and velocity) and resultant linear models are used in the controller designs [Topland (2004), Antonsen (2004), Blindheim (2004), Doruk (2008)]. An example solution by the linear quadratic regulator (LQR) is given by [Doruk (2008)] for the models with Modified Rodriguez Parameters. Linear quadratic methods provide easy and fast control design approaches and those properties can be applied to certain forms of nonlinear systems. The state dependent control parametrization is proposed for this purpose [Cloutier (1996), Cimen (2008)]. An application of the SDRE based approach to the satellite or spacecraft attitude control is presented by [Parrish (1997), Luo (2002)] where the attitude kinematics are modeled by the quaternion. All those solutions are suboptimal since the optimality conditions derived by [Cloutier (1996)] are difficult to be validated for high order nonlinear models where the satellite is one member. In this work, a different suboptimal SDRE based controller is proposed for the satellites where the kinematics is modeled by the Modified Rodriguez Parameters [Wiener (1962)].

# Satellite Model

In order to model the satellite one has to define the coordinate frames and their notations that are used throughout this paper. There are basically three coordinate frames that the satellite is referenced during its operation. Those are:

## *Inertial Reference Frame (Earth Centered Inertial or ECI):*

This is the primary reference coordinate system whose origin is fixed at the center of the earth. The z – axis coincides the earth's axis of rotation and points towards the north celestial pole. The x axis directed towards the vernal equinox and y axis completes the right hand rule. It is denoted by *i* and the components are thus shown by $x_i, y_i \& z_i$.

## *Orbit Referenced Frame:*

This frame has its origin on the mass center of the satellite. Its z – axis points towards the center of the earth and x – axis points in the travelling direction of the satellite. Lastly, the y – axis completes the right hand rule. Satellite rotations around those axes are sometimes called as roll, pitch and yaw rotations respectively. This frame is denoted by *o* and its components are $x_o, y_o \& z_o$.

## *Body Fixed Reference Frame:*

This coordinate system is fixed on the satellite's body and generally selected to coincide with the principal axes of inertia. Its components are denoted as $x_b, y_b \& z_b$ respectively.

# The attitude kinematics

The modified Rodriguez parameters are the minimal representations of attitude and thus have a singularity for the rotations of $\pm 360°$. This limitation does not produce a problem for practical ranges of rotations. The MRP is a three dimensional representation of the form:

$$\boldsymbol{\sigma} = \begin{bmatrix} \sigma_1 & \sigma_2 & \sigma_3 \end{bmatrix}^T \in \mathbb{R}^3 \quad (1)$$

The kinematics of this representation is:

$$\dot{\boldsymbol{\sigma}} = \mathbf{G}(\boldsymbol{\sigma})\boldsymbol{\omega} \quad (2)$$

The above kinematics represents the orientation of a three dimensional object with respect to a specific coordinate frame where the angular velocity vector $\boldsymbol{\omega} \in \mathbb{R}^3$ is also defined relative to the same coordinate frame. The matrix $\mathbf{G}(\boldsymbol{\sigma}) \in \mathbb{R}^{3\times 3}$ is always invertible and defined as shown below:

$$\mathbf{G}(\boldsymbol{\sigma}) = \frac{1}{2}\left( \frac{1-\boldsymbol{\sigma}^T\boldsymbol{\sigma}}{2} \mathbf{I}_{3\times 3} + \mathbf{S}(\boldsymbol{\sigma}) + \boldsymbol{\sigma}\boldsymbol{\sigma}^T \right) \quad (3)$$

Where the notation $\mathbf{S}(\boldsymbol{\sigma}) \in \mathbb{R}^{3\times 3}$ is the skew – symmetric operator:

$$\mathbf{S}(\mathbf{v}) = \begin{bmatrix} 0 & -v_3 & v_2 \\ v_3 & 0 & -v_1 \\ -v_2 & v_1 & 0 \end{bmatrix} \quad (4)$$

Another important parameter in the kinematics is the rotation matrix which defines a rotation from a specific frame to another specific frame. Mathematically, the rotation is expressed as a matrix multiplication and defined as $\mathbf{B}_2 = \mathbf{R}\mathbf{B}_1$ where $\mathbf{R} \in \mathbb{R}^{3\times3}$ is the rotation matrix from frame 1 to 2. $\mathbf{B}_1 \& \mathbf{B}_2$ are any vectorial measurements defined in the frames 1 and 2 respectively. In terms of the modified Rodriguez parameters the rotational operation generally takes the following form:

$$\mathbf{B}_B = \mathbf{R}(\sigma)\mathbf{B}_I \quad (5)$$

where $\mathbf{R}(\sigma) \in \mathbb{R}^{3\times3}$ is the rotation matrix from the inertial to the body frame, $B_I \in \mathbb{R}^3$ is any reference vector defined on the inertial reference frame and the $B_B \in \mathbb{R}^3$ is another reference vector defined in the satellite body frame. The mathematical expression of $\mathbf{R}(\sigma)$ is given below:

$$\mathbf{R}(\sigma) = \mathbf{I}_{3\times3} - \frac{4(1-\sigma^T\sigma)}{(1+\sigma^T\sigma)^2}\mathbf{S}(\sigma) + \frac{8}{(1+\sigma^T\sigma)^2}\mathbf{S}^2(\sigma) \quad (6)$$

In the satellite attitude control approaches the satellite attitude can be defined relative to the orbit frame. In that case the kinematics equation is rewritten as:

$$\dot{\sigma}_e = \mathbf{G}(\sigma_e)\omega_{ob}^b \quad (7)$$

In the above $\sigma_e$ is the attitude of the satellite relative to the orbital frame and $\omega_{ob}^b$ is the angular velocity of the satellite body relative to the orbital frame. The latter is defined mathematically as:

$$\omega_{ob}^b = \omega_{ib}^b - \mathbf{R}(\sigma_e)\omega_{io}^o \quad (8)$$

Where $\omega_{ib}^b$ is the angular velocity of the satellite body relative to the inertial reference frame (ECI) and $\omega_{io}^o$ is the orbital velocity of the satellite. The rotation matrix $\mathbf{R}(\sigma_e)$ is used to transform the orbital velocity to the satellite body frame. If the orbital frame is considered as the desired reference frame this will lead to a tracking error problem since the relative attitude of the satellite reduces to the attitude error with respect to the orbital reference frame. So one can denote the angular velocity term $\omega_{ob}^b$ as $\omega_e$ in order to stress that it is an angular velocity error. Finally, the kinematics equation of (7) can be rewritten once as:

$$\dot{\sigma}_e = \mathbf{G}(\sigma_e)\omega_e \quad (9)$$

### Satellite dynamics:

For a rigid body satellite, the dynamics can be derived from the Newton's laws [Fauske (2002)]. The resultant model in terms of the body angular velocity is a nonlinear differential equation shown as:

$$\dot{\boldsymbol{\omega}}_{ib}^{b} = -\mathbf{I}^{-1}\mathbf{S}\left(\boldsymbol{\omega}_{ib}^{b}\right)\mathbf{I}\boldsymbol{\omega}_{ib}^{b} + \mathbf{I}^{-1}\boldsymbol{\tau}_{a} \tag{10}$$

where $\boldsymbol{\tau}_a \in \mathbb{R}^3$ is the actuator torque input and $\mathbf{I} \in \mathbb{R}^{3\times 3}$ is the symmetric moment of inertia matrix which can be shown as:

$$\mathbf{I} = \begin{bmatrix} I_x & I_{xy} & I_{xz} \\ I_{xy} & I_y & I_{yz} \\ I_{xz} & I_{yz} & I_z \end{bmatrix} \tag{11}$$

If one desires to represent the dynamics in terms of the error angular velocity $\boldsymbol{\omega}_e$, the resultant equation is:

$$\begin{aligned}\dot{\boldsymbol{\omega}}_e &= -\mathbf{I}^{-1}\mathbf{S}\left(\boldsymbol{\omega}_e + \mathbf{R}(\boldsymbol{\sigma}_e)\boldsymbol{\omega}_{io}^o\right)\mathbf{I}\left[\boldsymbol{\omega}_e + \mathbf{R}(\boldsymbol{\sigma}_e)\boldsymbol{\omega}_{io}^o\right] - \dot{\mathbf{R}}(\boldsymbol{\sigma}_e)\boldsymbol{\omega}_{io}^o + \mathbf{I}^{-1}\boldsymbol{\tau}_a \\ \dot{\mathbf{R}}(\boldsymbol{\sigma}_e) &= \mathbf{S}(\boldsymbol{\omega}_e)\mathbf{R}(\boldsymbol{\sigma}_e)\end{aligned} \tag{12}$$

When the above equation is written in its expanded form:

$$\dot{\boldsymbol{\omega}}_e = -\mathbf{I}^{-1}\mathbf{S}(\boldsymbol{\omega}_e)\mathbf{I}\boldsymbol{\omega}_e - \mathbf{I}^{-1}\mathbf{S}(\mathbf{R}(\boldsymbol{\sigma}_e)\boldsymbol{\omega}_{io}^o)\mathbf{I}\boldsymbol{\omega}_e - \mathbf{I}^{-1}\mathbf{S}(\boldsymbol{\omega}_e)\mathbf{I}\mathbf{R}(\boldsymbol{\sigma}_e)\boldsymbol{\omega}_{io}^o - \mathbf{I}^{-1}\mathbf{S}(\mathbf{R}(\boldsymbol{\sigma}_e)\boldsymbol{\omega}_{io}^o)\mathbf{I}\mathbf{R}(\boldsymbol{\sigma}_e)\boldsymbol{\omega}_{io}^o - \mathbf{S}(\boldsymbol{\omega}_e)\mathbf{R}(\boldsymbol{\sigma}_e)\boldsymbol{\omega}_{io}^o + \mathbf{I}^{-1}\boldsymbol{\tau}_a \tag{13}$$

For the sake of simplicity one can write the above in the following form:

$$\dot{\boldsymbol{\omega}}_e = -\mathbf{I}^{-1}\mathbf{S}(\boldsymbol{\omega}_e)\mathbf{I}\boldsymbol{\omega}_e + \mathbf{I}^{-1}\mathbf{u} \tag{14}$$

where,

$$\mathbf{u} = -\mathbf{S}(\mathbf{R}(\boldsymbol{\sigma}_e)\boldsymbol{\omega}_{io}^o)\mathbf{I}\boldsymbol{\omega}_e - \mathbf{I}\mathbf{R}(\boldsymbol{\sigma}_e)\boldsymbol{\omega}_{io}^o - \mathbf{S}(\boldsymbol{\omega}_e)\mathbf{I}\mathbf{R}(\boldsymbol{\sigma}_e)\boldsymbol{\omega}_{io}^o - \mathbf{S}(\mathbf{R}(\boldsymbol{\sigma}_e)\boldsymbol{\omega}_{io}^o)\mathbf{I}\mathbf{R}(\boldsymbol{\sigma}_e)\boldsymbol{\omega}_{io}^o - \mathbf{I}\mathbf{S}(\boldsymbol{\omega}_e)\mathbf{R}(\boldsymbol{\sigma}_e)\boldsymbol{\omega}_{io}^o + \boldsymbol{\tau}_a \tag{15}$$

$\mathbf{u}$ can be considered as a virtual input where the control law will be derived on. In a possible practical application one should extract $\boldsymbol{\tau}_a$ from the control law.

In the next section, the state dependent Riccati equation based controller design technique is to be introduced.

## State Dependent Riccati Equation (SDRE) Technique

Consider the nonlinear system shown below:

$$\dot{\mathbf{x}} = \mathbf{f}(\mathbf{x},\mathbf{u}),\ \mathbf{x} \in \mathbb{R}^n,\ \mathbf{u} \in \mathbb{R}^m \tag{16}$$

In some cases the governing system function $f(x,u)$ can be written as $f(x,u) = A(x)x + B(x)u$ where $A(x) \in \mathbb{R}^{n \times n}$ and $B(x) \in \mathbb{R}^{n \times m}$. So the new system model is:

$$\dot{x} = A(x)x + B(x)u \quad (17)$$

If one desires to minimize the state dependent quadratic cost function:

$$J = \int_0^\infty \{x^T Q(x) x + u^T R(x) u\} dt \quad (18)$$

The solution of the following equation will solve the problem by providing a control law on $u$ as $u = -K(x)x$ with $K(x) = -R^{-1}(x) B^T(x) P(x) x$:

$$A^T(x) P(x) + P(x) A(x) - B(x) P(x) R^{-1}(x) P(x) B^T(x) = -Q(x) \quad (19)$$

where $P(x) \in \mathbb{R}^{n \times n}$ is the positive definite symmetric solution of the above equation. The matrices $Q(x)$ and $R(x)$ are also positive definite and symmetric. In some cases, they can be taken as constants. If a positive definite and symmetric $P(x) \in \mathbb{R}^{n \times n}$ is found then the closed loop formed by using the control law $u = -K(x)x$ is point wisely stable at each operating point $x = x_o$.

## SDRE based Attitude Control:

Before going into the main issue, it will be convenient to repeat the satellite attitude kinematics and dynamics for convenience:

$$\begin{aligned} \dot{\sigma}_e &= G(\sigma_e) \omega_e \\ \dot{\omega}_e &= -I^{-1} S(\omega_e) I \omega_e + I^{-1} u \end{aligned} \quad (20)$$

For usage in SDRE technique one should convert the above model into the form given in (17):

$$\begin{bmatrix} \dot{\sigma}_e \\ \dot{\omega}_e \end{bmatrix} = \begin{bmatrix} 0 & G(\sigma_e) \\ 0 & -I^{-1} S(\omega_e) I \end{bmatrix} \begin{bmatrix} \sigma_e \\ \omega_e \end{bmatrix} + \begin{bmatrix} 0 \\ I^{-1} \end{bmatrix} u \quad (21)$$

So:

$$A(x) = \begin{bmatrix} 0 & G(\sigma_e) \\ 0 & -I^{-1} S(\omega_e) I \end{bmatrix}, B(x) = \begin{bmatrix} 0 \\ I^{-1} \end{bmatrix}, x = \begin{bmatrix} \sigma_e \\ \omega_e \end{bmatrix} \quad (22)$$

The matrix $B(x)$ is not state dependent in the current case. The equation in (19) is solved for the following operating points:

$$\boldsymbol{\sigma}_e \to [-1, \ 0, \ 1]$$
$$\boldsymbol{\omega}_e \to [-5, \ 0, \ 5]°/\sec \tag{23}$$

After the computation each component of the control gain matrix $\mathbf{K}(\mathbf{x})$ is put in a lookup table of dimension equal to six. During the operation they are interpolated linearly (or in another way chosen by the designer). If the control gain matrix is denoted as shown below:

$$\mathbf{K}(\mathbf{x}) = \begin{bmatrix} k_{11}(\mathbf{x}) & k_{12}(\mathbf{x}) & k_{13}(\mathbf{x}) & k_{14}(\mathbf{x}) & k_{15}(\mathbf{x}) & k_{16}(\mathbf{x}) \\ k_{21}(\mathbf{x}) & k_{22}(\mathbf{x}) & k_{23}(\mathbf{x}) & k_{24}(\mathbf{x}) & k_{25}(\mathbf{x}) & k_{26}(\mathbf{x}) \\ k_{31}(\mathbf{x}) & k_{31}(\mathbf{x}) & k_{31}(\mathbf{x}) & k_{34}(\mathbf{x}) & k_{35}(\mathbf{x}) & k_{36}(\mathbf{x}) \end{bmatrix} \tag{24}$$

each $k_{ij}(\mathbf{x})$ can be expressed as a lookup table having the indices as the attitude and angular velocity errors. The final control law is:

$$\begin{aligned}\boldsymbol{\tau}_a = &-\mathbf{K}(\mathbf{x})\mathbf{x} + \mathbf{S}(\mathbf{R}(\boldsymbol{\sigma}_e)\boldsymbol{\omega}_{io}^o)\mathbf{I}\boldsymbol{\omega}_e + \mathbf{S}(\boldsymbol{\omega}_e)\mathbf{I}\mathbf{R}(\boldsymbol{\sigma}_e)\boldsymbol{\omega}_{io}^o + \mathbf{S}(\mathbf{R}(\boldsymbol{\sigma}_e)\boldsymbol{\omega}_{io}^o)\mathbf{I}\mathbf{R}(\boldsymbol{\sigma}_e)\boldsymbol{\omega}_{io}^o \\ &+ \mathbf{I}\mathbf{S}(\boldsymbol{\omega}_e)\dot{\mathbf{R}}(\boldsymbol{\sigma}_e)\dot{\boldsymbol{\omega}}_{io}^o \end{aligned} \tag{25}$$

## Simulation Example

### Outline of the Simulation:
The model parameters of the satellite are taken from a realistic satellite application [Kaplan (2006)] where the inertia matrix is the following:

$$\mathbf{I} = \begin{bmatrix} 9.8194 & -0.0721 & -0.2893 \\ -0.0721 & 9.7030 & -0.1011 \\ -0.2893 & -0.1011 & 9.7309 \end{bmatrix} kg \cdot m^2$$

The satellite starts from the attitude $\boldsymbol{\sigma}_{initial} = [0.3333 \ -0.3333 \ -0.3333]^T$ and moves to the zero attitude position. In Euler angles that corresponds to the trajectory starting from $\boldsymbol{\phi} = [90° \ 0° \ -90°]^T$ and ending at the origin. Two cases of simulations are presented in this section. In the first case, the simulation is performed without the existence of any parametric uncertainty whereas in the second case there are parametric uncertainties on the inertia matrix. In the uncertainty existent environment the simulations are repeated a hundred times in order to achieve reliable information.

### The quadratic performance indices:
The quadratic performance indices are taken as constant matrices in this research. The forms of the matrices are given below for convenience:

$$\mathbf{Q}(\mathbf{x}) = \mathbf{Q} = \begin{bmatrix} q_1 \mathbf{I}_{3\times 3} & 0 \\ 0 & q_2 \mathbf{I}_{3\times 3} \end{bmatrix} \quad (26)$$

$$\mathbf{R}(\mathbf{x}) = \mathbf{R} = r\mathbf{I}_{3\times 3}$$

For the current simulation four cases are selected which are shown below:

Table 1 The values of quadratic weighting coefficients and comments on the simulation results

| No. of Case | $q_1$ | $q_2$ | $r$ | Comments |
|---|---|---|---|---|
| 1 | 0.1 | 0.000001 | 1 | $\max(\boldsymbol{\tau}_a) \geq 0.02 N \cdot m$ |
| 2 | 0.1 | 0.000001 | 1000 | $\boldsymbol{\tau}_a$ normal |
| 3 | 0.001 | 0.000001 | 1 | $\boldsymbol{\tau}_a$ normal |
| 4 | 0.001 | 0.000001 | 1000 | $\boldsymbol{\tau}_a$ low but the system response is too slow |

In the table above, some comments on the simulation results are also given, the comments on the torque levels are given according to the characteristics of the satellite model presented in [Kaplan (2006)] where it is advised not to have a torque requirement larger than 0.02 Nm. However, there is no actuator limit set on the actuator torque control law. In the next section the graphical results are given.

## Graphical Results:

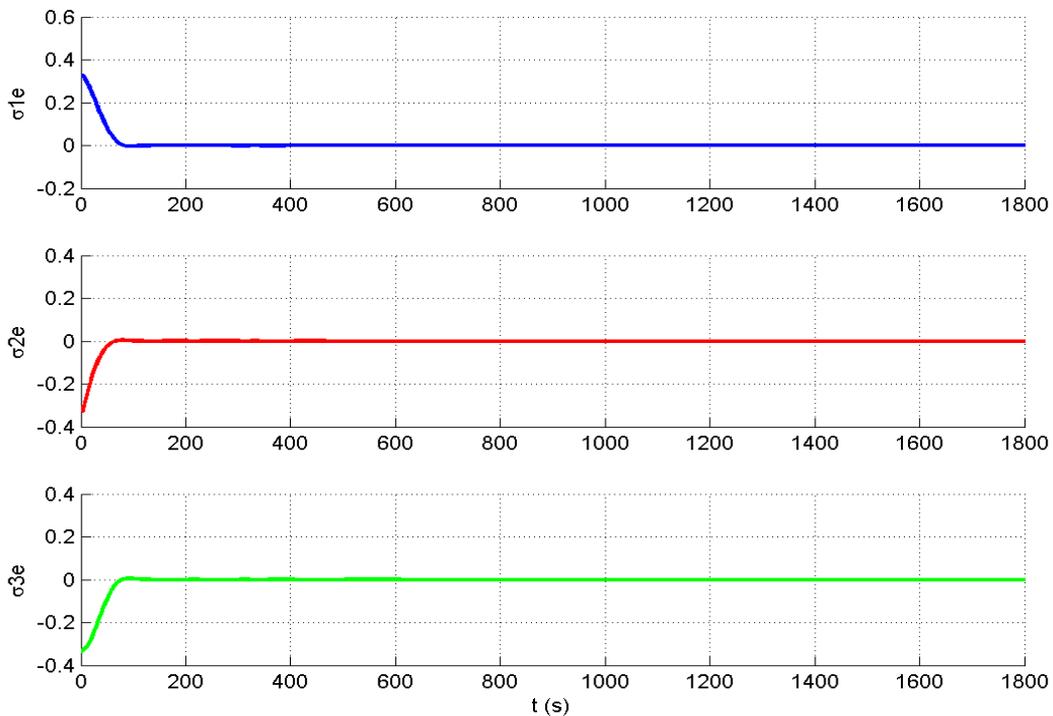

Figure 1 Attitude error for the first case

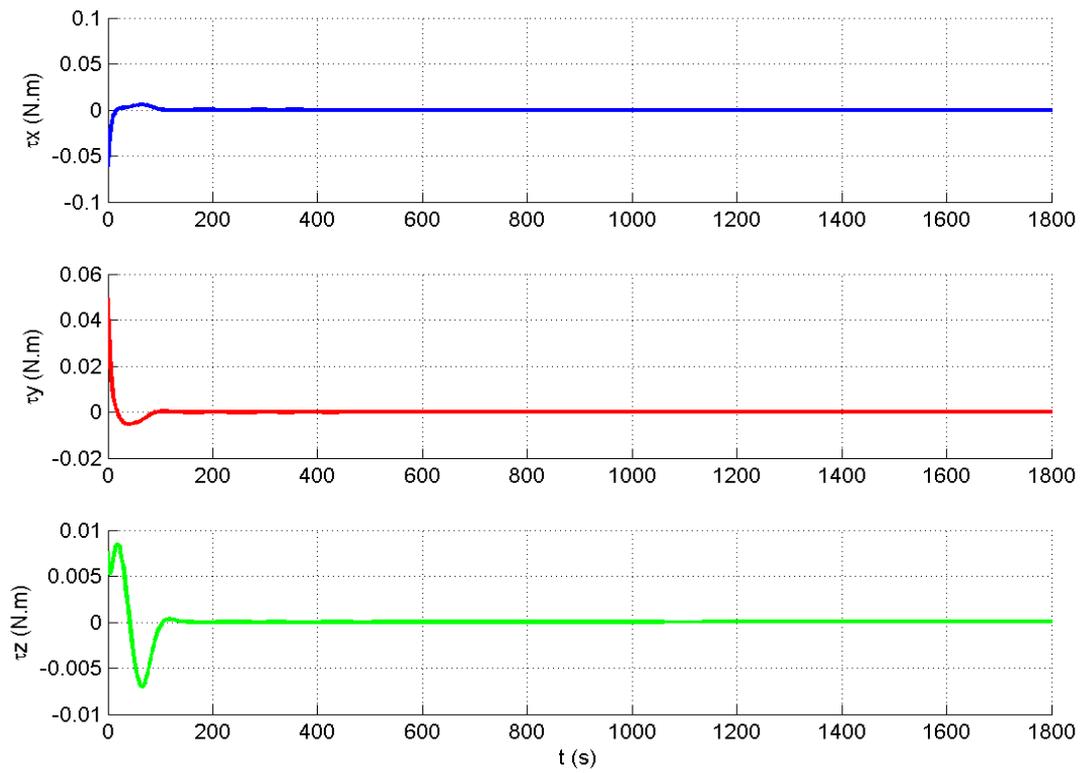

**Figure 2 Torque requirements for the first case**

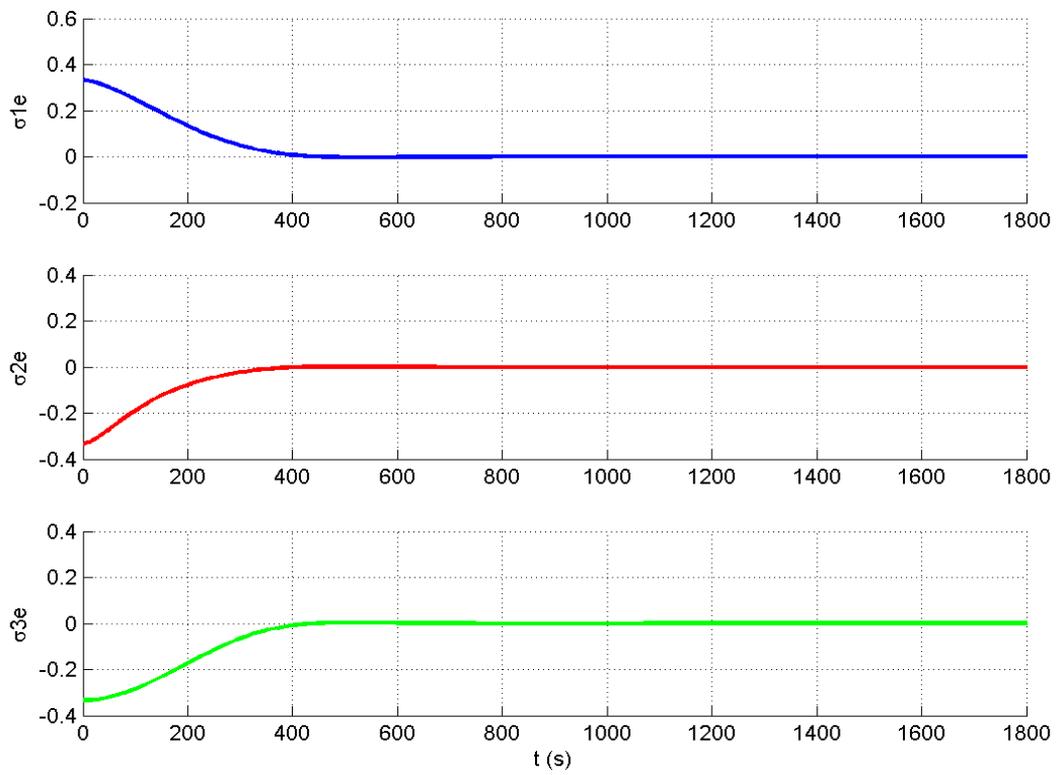

**Figure 3 Attitude error for the second case**

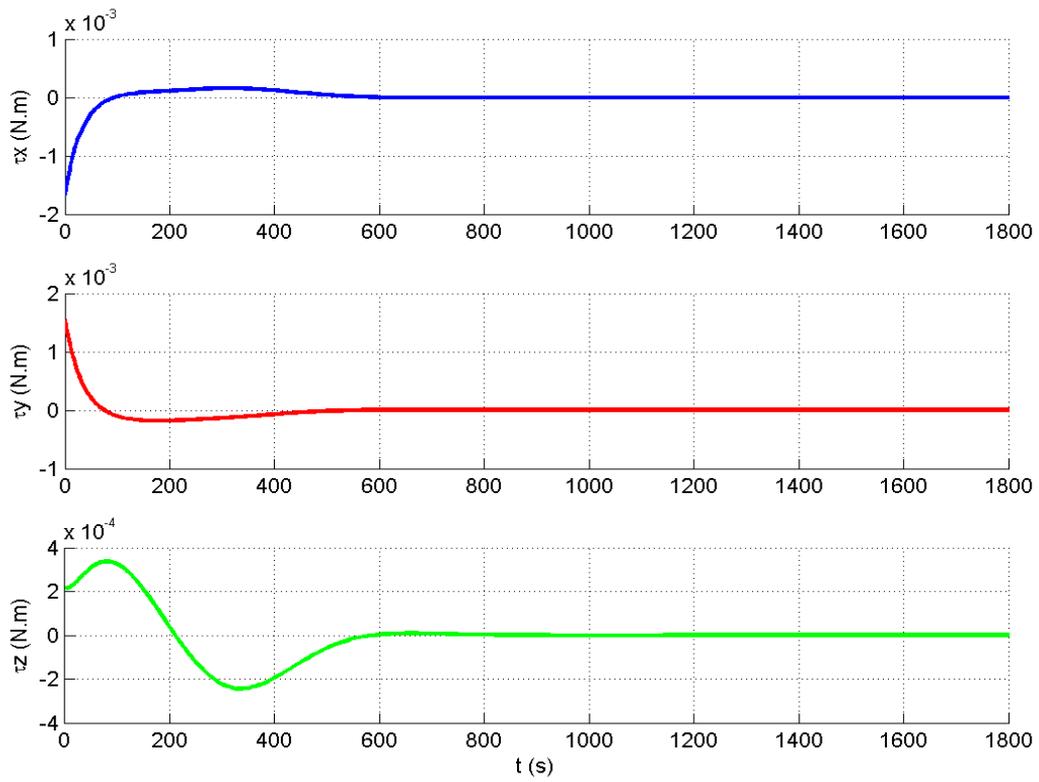

Figure 4 Torque requirements for the second case

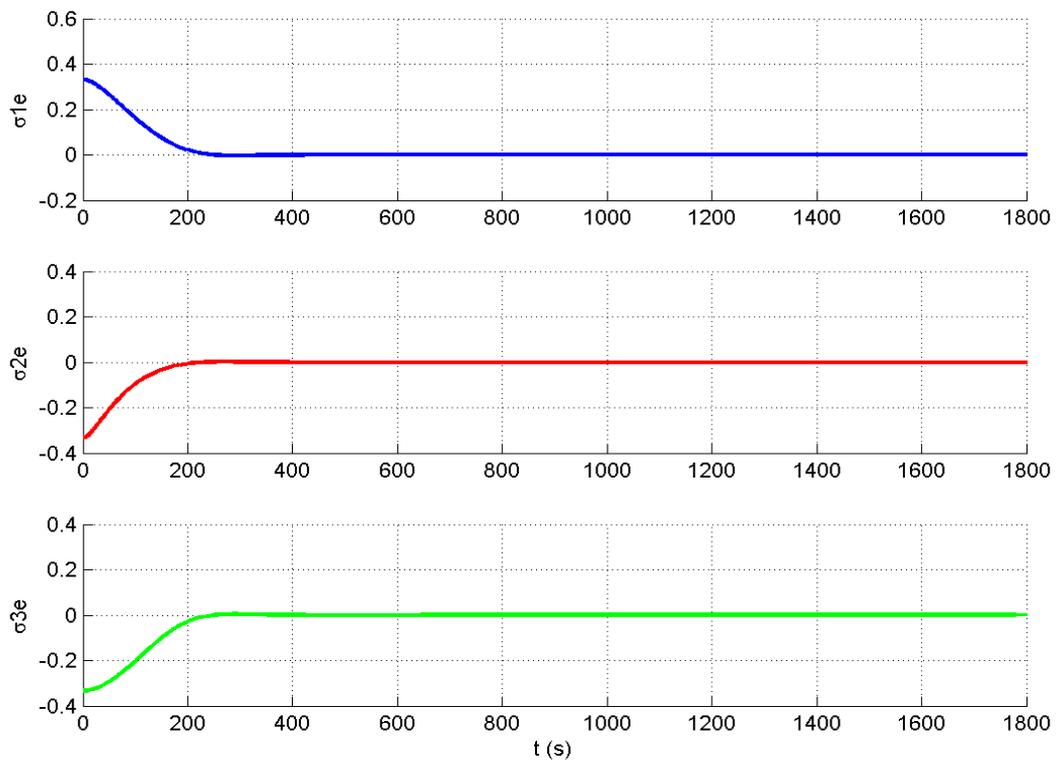

Figure 5 Attitude error for the third case

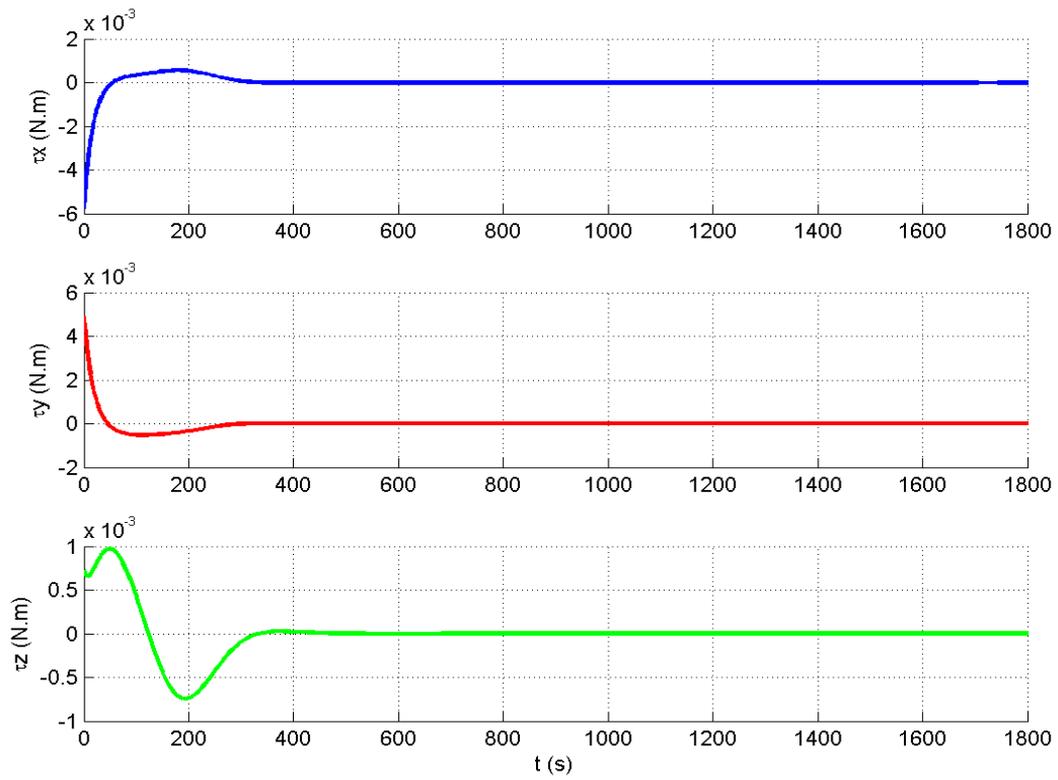

Figure 6 Torque requirement for the third case

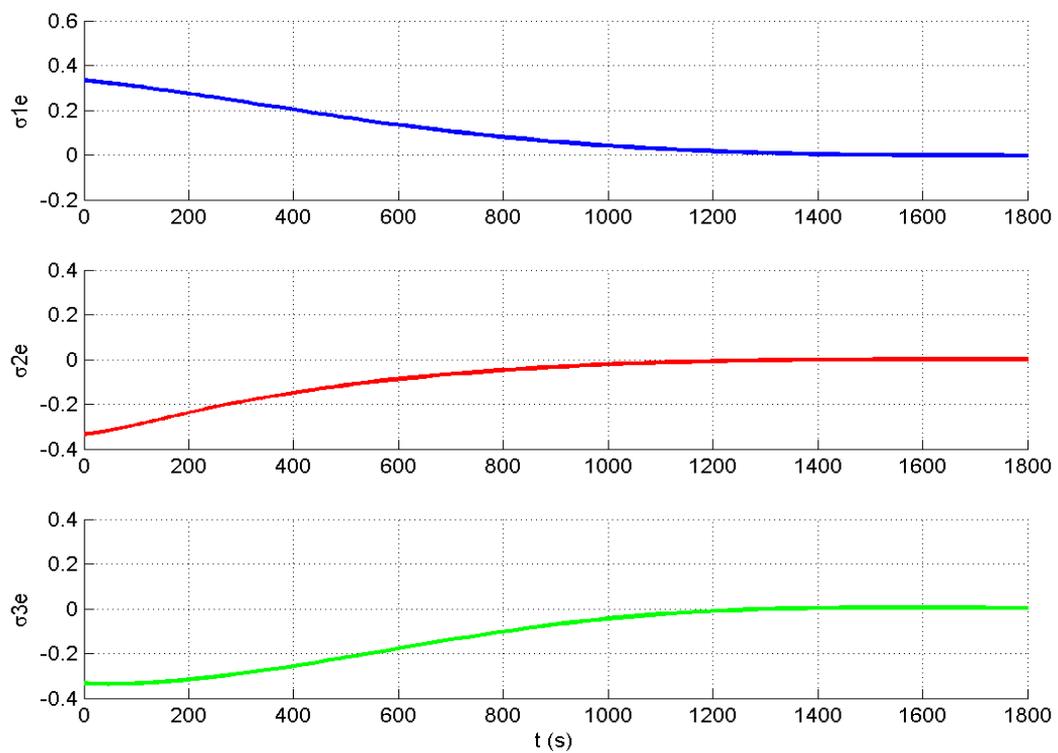

Figure 7 Attitude error for the fourth case

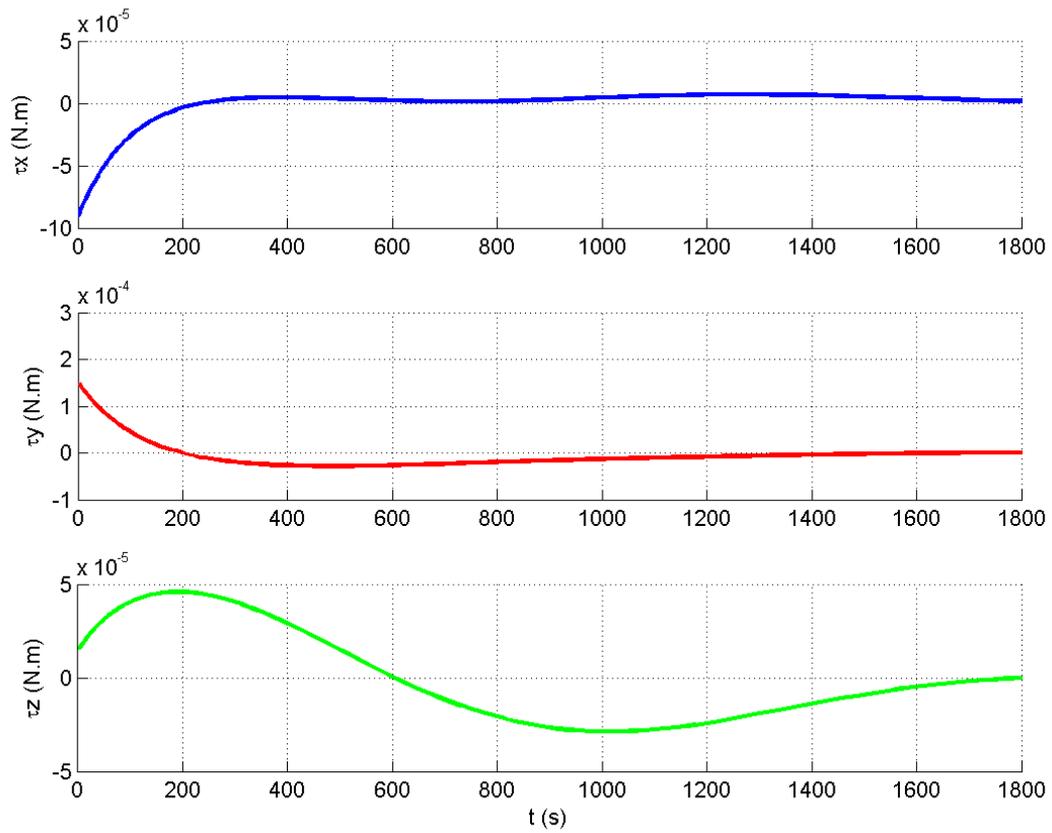

**Figure 8 Torque requirement for the fourth case**

As it can be understood from the graphical results the best candidate for application seems to be the second entry in Table 1. The torque requirement of the fourth case seems lower but its response is too slow and the torque levels of the second case is far under 0.02 Nm so it will be the best to use the second weight set in Table 1. Below the results of the simulation in presence of parametric uncertainties on the inertia matrices:

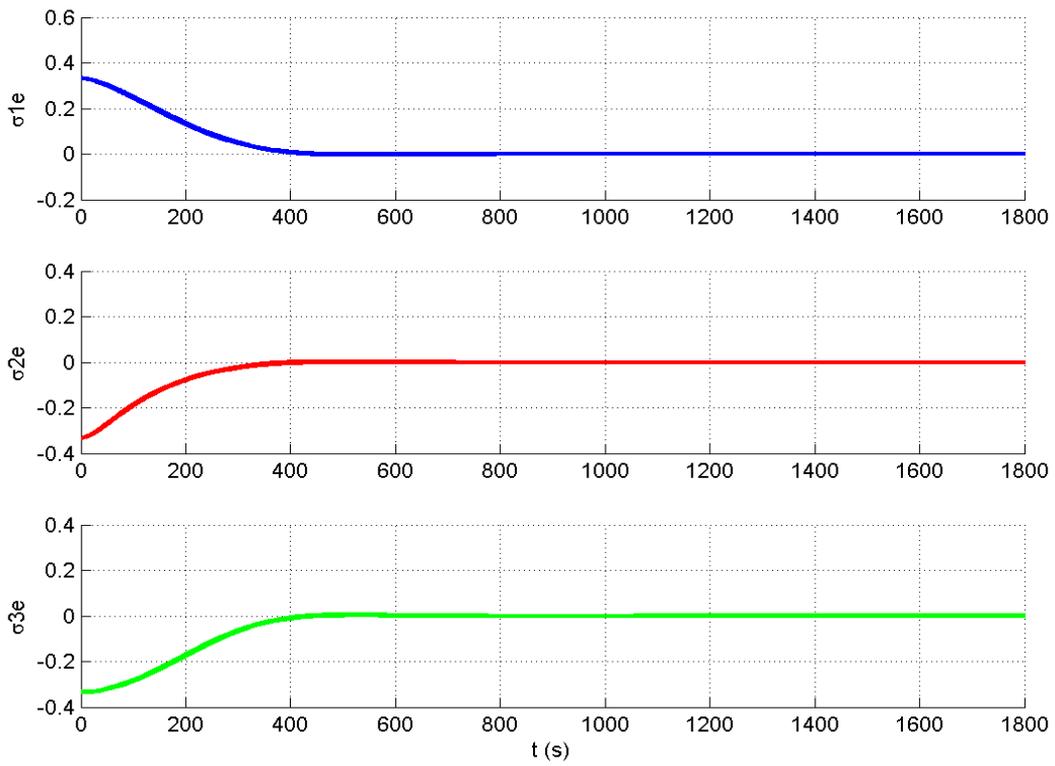

**Figure 9 Attitude error for the second case in presence of inertia uncertainties**

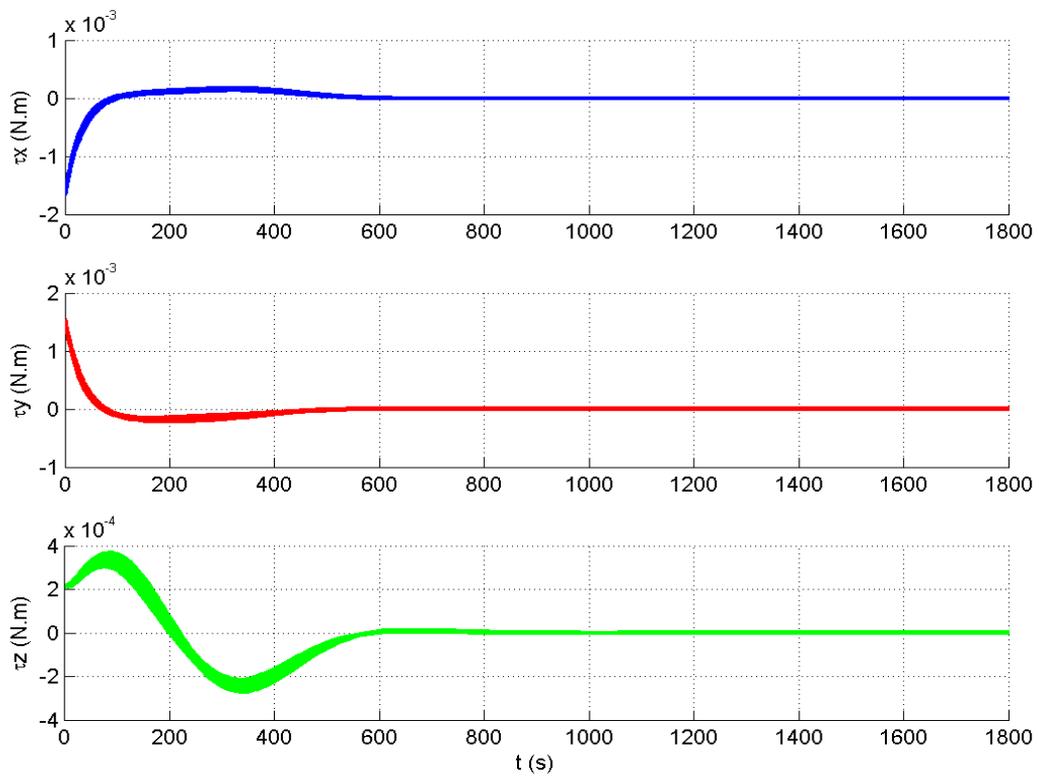

**Figure 10 Torque requirements for the second case in presence of inertia uncertainties.**

# Conclusion:

In this research we have designed an optimal attitude controller based on state dependent Riccati equation (SDRE) technique. The satellite attitude kinematics is represented by the Modified Rodriguez Parameters (MRP) which is minimalist in nature. Due to the difficulty of finding an analytical solution for the control law certain points are taken in the operating range and a lookup table is created. The necessary control gains are evaluated in the selected points and are interpolated throughout the operation of the satellite. The positive definite solution of the SDRE is found by the aid of MATLAB Control System Toolbox and the resultant control gains form point wisely stable control laws. The simulations results show that the attitude control laws operate successfully even in the presence of inertia matrix uncertainties. Analysis of global stability of the closed loop controller may constitute the subject of possible future research.